\documentclass[12pt,a4paper,reqno]{amsart}

\usepackage{color}
\usepackage{enumerate}
\usepackage[margin=2.5cm]{geometry}
\usepackage{graphicx} 
\usepackage{overpic}

\usepackage{hyperref}

\newtheorem{lema}{Lemma}[section]
\newtheorem{tho}[lema]{Theorem}
\newtheorem{pro}[lema]{Proposition}

\newtheorem{coro}[lema]{Corollary}

\DeclareMathOperator{\Jac}{Jac}

\newcommand{\R}{\ensuremath{\mathbb{R}}}

\begin{document}
	

\title[Limit cycles in Kolmogorov quadratic piecewise systems]{Limit cycles in piecewise \\ quadratic Kolmogorov systems}

\author[L.P.C. da Cruz]
{Leonardo Pereira Costa da Cruz}
\address{Instituto de Ciências Matemáticas e de Computaç\~o. Universidade de S\~o Paulo, S\~ao Carlos, 13566--590 S\~ao Paulo, Brazil}
\email{leonardocruz@icmc.usp.br}

\author[R. Oliveira]
{Regilene Oliveira}
\address{Instituto de Ciências Matemáticas e de Computaç\~o. Universidade de S\~o Paulo, S\~ao Carlos, 13566--590 S\~ao Paulo, Brazil}
\email{regilene@icmc.usp.br}

\author[J. Torregrosa]{Joan Torregrosa}
\address {Departament de Matem\`{a}tiques, Facultat de Ci\`{e}ncies, Universitat Aut\`{o}noma de Bar\-ce\-lo\-na, 08193 Bellaterra, Barcelona, Spain ; and Centre de Recerca Matem\`{a}tica, Edifici Cc, Campus de Bellaterra, 08193 Cerdanyola del Vall\`{e}s (Barcelona), Spain}
\email{joan.torregrosa@uab.cat}

	

\begin{abstract}
We study planar piecewise quadratic differential systems of Kolmogorov type. Specifically, we consider systems with both coordinate axes invariant and with a separation line that is straight and distinct from the invariant axes. The main results concern two different aspects. First, the center problem is solved for certain subclasses. Second, using this classification, the bifurcation of limit cycles of crossing type is investigated. We contrast the nature of Hopf-type bifurcations in smooth and piecewise smooth settings, particularly regarding the bifurcation of limit cycles of small amplitude. The classical pseudo-Hopf bifurcation is analyzed in the Kolmogorov systems class. It is worth highlighting that, in contrast to the smooth Kolmogorov quadratic systems, which have no limit cycles, the piecewise case exhibits at least six. Furthermore, we show that the maximal weak focus order, eight, does not necessarily yield the maximal number of small-amplitude limit cycles.
\end{abstract}
	
\maketitle
	
\section{Introduction}
There has been significant interest in the study of piecewise differential systems in recent years. This interest is likely due to the fact that many natural phenomena can be modeled by such systems. Examples include electrical and mechanical systems, control theory, and even genetic networks; see \cite{AcaBonBro2011,BerBudCha2008,Fil1988}. Additionally, from a theoretical point of view, many authors are now exploring the piecewise setting by posing questions analogous to those studied in the smooth case, such as determining an analogue of the Hilbert number for piecewise polynomial systems, characterizing centers, investigating conditions for integrability, and so on.

We adopt this second perspective, focusing on the extension of the concepts of cyclicity and the center problem within the class of planar piecewise differential systems of the form
\[
Z = (Z_1(x,y), Z_2(x,y)),
\]
where $Z_i$, for $i = 1, 2$, are smooth vector fields defined respectively in the regions $\Sigma_i = \{(x,y) : (-1)^i h(x,y) > 0\} $, with $h: \mathbb{R}^2 \to \mathbb{R}$ a $\mathcal{C}^1$ function for which 0 is a regular value. In this setting, the curve $\Sigma = \{ h(x,y) = 0 \} $ is called the separation line. More specifically, we investigate the behavior of closed solutions of $Z$ near pseudo-equilibrium points with monodromic behavior, i.e. points where solutions of $Z$ rotate around the point. In this context, two possibilities arise: the study of isolated periodic orbits (limit cycles) and of continua of periodic orbits (centers).

It is worth noting that the dynamics of polynomial piecewise differential systems of degree $n$ are richer than those of smooth systems of the same degree. For instance, it is well known that linear differential systems do not admit limit cycles. However, in piecewise linear systems with two zones separated by a straight line, examples exist with at least three limit cycles. A configuration with three limit cycles was first detected numerically by Huan and Yang in \cite{HuanYang2010}, and later confirmed analytically by Llibre and Ponce in \cite{LliPon2012}. The existence of three limit cycles has also been obtained via perturbations of a center: Buzzi et al. \cite{BuzPesTor2013} found three limit cycles arising from a fourth-order piecewise linear perturbation of a linear center, while Llibre et al. \cite{LliNovTei2015a} achieved the same result using a first-order perturbation of a specific piecewise linear center. See also Freire et al. \cite{FrePonTor2014} for bifurcating them from infinity.

Bautin \cite{Bautin1954} proved that planar quadratic systems can have at most three small-amplitude limit cycles bifurcating from a monodromic equilibrium, and this remains the only complete classification, with respect to degree, for which the center-focus problem is fully solved. Recently, at least twelve limit cycles have been found for polynomial piecewise quadratic systems; see \cite{BraCruTor2024}.

A similar situation occurs in quadratic differential systems with an invariant straight line. Cherkas, Zhilevich, and Rychkov \cite{CheZil1970,CheZil1972,Ryc1972} proved that such systems have at most one limit cycle. This problem was revisited in \cite{ColLli1990}. In the corresponding piecewise class, however, it was shown in \cite{daCruzTor2022} that systems with at least seven limit cycles exist. Recall also that quadratic systems with two transversal invariant straight lines have no limit cycles but rather a continuum of periodic orbits, as in the Lotka--Volterra class. It is well known that the absence of limit cycles in this case is due to the existence of a Dulac function when the system is not integrable. This result can be found in several classical textbooks, such as \cite{Chicone2024, Perko2001}. 

In this paper, we focus on planar piecewise quadratic differential systems with two transversal invariant straight lines, i.e. a natural restriction on the general class of piecewise autonomous planar Kolmogorov systems. As usual, after an affine change of coordinates if necessary, we can assume that these systems can be written as
\[
(\dot{x}, \dot{y}) = (xP(x,y),\ yQ(x,y)),
\]
where $P$ and $Q$ are polynomials of degree one in the variables $x$ and $y$. Recently, a particular case which exhibits at least one limit cycle and known as Palomba's model has beeen studied in \cite{YagGouOleg2025}. Kolmogorov systems are frequently used to model the interaction of two species occupying the same ecological niche \cite{CoeRebSov2021, Zhao2017}. In fact, many natural phenomena can be modeled by Kolmogorov systems, for instance, in mathematical ecology and population dynamics, chemical reactions, plasma physics, hydrodynamics, etc. Probably the most well-known Kolmogorov systems are the Lotka--Volterra systems \cite{Lotka1925, Volterra1927}, which correspond to the case where $P(x,y) = \alpha - \beta y$ and $Q(x,y) = -\gamma + \delta x$, with all parameters positive. 

Although there are many works on the dynamics of Kolmogorov systems, as far as the authors know, there are no results concerning the bifurcation of limit cycles in piecewise quadratic Kolmogorov systems.

In this work, we prove the existence of piecewise Kolmogorov quadratic systems separated by a straight line that have at least six limit cycles (see Theorem \ref{th:m1}). More concretely, we consider planar piecewise systems in the set $\mathcal{K}_2$, the class of Kolmogorov quadratic differential systems defined in two zones separated by a straight line, i.e., systems of the form
\begin{equation}\label{eq:2}
	Z_i:\begin{cases}
		\dot{x} = x (a_i + b_i x + c_i y),\\
		\dot{y} = y (d_i + e_i x + f_i y),
	\end{cases} \quad \text{if} \quad (x,y) \in \Sigma_i = {(x,y) : (-1)^i h(x,y) > 0}
\end{equation}
where $i = 1, 2$, and $\Sigma = \{ h(x,y) = 0 \}$, with $h$ a degree one polynomial. Observe that each system in $\mathcal{K}_2$ admits the two invariant straight lines $x = 0$ and $y = 0$.

Summarizing, the contributions of this paper lie in the study of an extension of the classical center-focus and cyclicity problems to non-smooth planar systems in the class $\mathcal{K}_2$. More specifically, we study the stability of monodromic equilibrium points on $\Sigma$, the order of weakness, and the number of crossing limit cycles that bifurcate from them. Recall that a crossing limit cycle is an isolated periodic orbit that intersects both zones defined by the separation line $\Sigma$ without a sliding segment. As in the smooth case, we say that such a limit cycle is of small amplitude if there exists a small neighborhood of the pseudo-equilibrium $p$ in which the limit cycle is entirely contained. More details on this bifurcation and the definition of pseudo-equilibrium will be given in what follows.

The main results of this paper are described below.

\begin{tho}\label{th:m1} There exist planar piecewise Kolmogorov quadratic differential systems, defined as in \eqref{eq:2}, having at least six crossing limit cycles of small amplitude.
\end{tho}

In Section~\ref{bifurcation}, we will see that the number of crossing limit cycles in $\mathcal{K}_2$ depends on the separation straight line $\Sigma$. Hence, we will provide different results depending on the choice of the separation line.

The definition of a weak focus order for an equilibrium point in the class of piecewise differential systems differs from the classical smooth scenario and will be introduced in Section~\ref{se:preliminaries}. In this section, we also discuss the notion of a degenerate Hopf bifurcation in the piecewise setting. In this context, even though we obtain limit cycles using a versal unfolding of a degenerate Hopf bifurcation point of sixth-order in Theorem~\ref{th:m1}, the problem of the highest weakness is more intricate. In fact, we prove the existence of weak focus points of eighth-order, but without a versal unfolding in $\mathcal{K}_2$, in this case, we obtain only five limit cycles.

The complete characterization of the center problem in the $\mathcal{K}_2$ class is very difficult; the next result gives a partial answer. We restrict our analysis to the case where both systems in \eqref{eq:2} have a common equilibrium point on $\Sigma$ with purely imaginary eigenvalues in their Jacobian matrices. It can be checked that these conditions are sufficient to guarantee that the equilibrium point of each system $Z_i$, $i = 1, 2$, is a center. We call such points of the $CC$-equilibrium type. In this context, $\det\Jac Z(x_0, y_0)$ denotes the determinant of the Jacobian matrix of $Z$ evaluated at the point $(x_0, y_0)$.

\begin{tho}\label{c3}
Let $Z \in \mathcal{K}_2$ be defined as in \eqref{eq:2}, and let $(x_0, y_0)$ be a monodromic $CC$-equilibrium point on $\Sigma$, with $\det \Jac Z_i(x_0, y_0) = D_i$ for $i = 1, 2$. Then, the next statements hold.
\begin{enumerate}[(i)]
\item If $h(x,y)=x y_0-y x_0$ then $Z$ has a center at $(x_0,y_0)$ if, and only if,  
$$D_1^2D_2^2(b_1e_2-b_2e_1)-x_0^2D_1^2b_2e_1(b_2 - e_2)^2 + x_0^2D_2^2b_1e_2(b_1 - e_1)^2=0.$$
\item If $h(x,y)=y-y_0$ then $Z$ has a center  at $(x_0,y_0)$  if, and only if, $$D_1^2b_2^2 - D_1^2b_2e_2 - D_2^2b_1^2 + D_2^2b_1e_1=0.$$
\item If $h(x,y)=x-x_0$ then $Z$ has a center at $(x_0,y_0)$  if, and only if, 			
$$D_1^2D_2^2(b_1e_2 - b_2e_1) - x_0^2D_1^2b_2^2e_1(b_2 - e_2) + x_0^2D_2^2b_1^2e_2(b_1 - e_1)=0.$$
\end{enumerate}
\end{tho}
This paper is structured as follows. In Section~\ref{se:preliminaries}, we present different types of bifurcations of small-amplitude limit cycles, namely the well-known Hopf and pseudo-Hopf types. The proof of Theorem~\ref{c3} and the characterization of centers for a specific class of $\mathcal K_2$ are given in Section~\ref{centers}. Finally, in Section~\ref{bifurcation}, the bifurcation of limit cycles in $\mathcal{K}_2$ is investigated, which includes the proof of Theorem~\ref{th:m1}.

\section{Preliminaries on Hopf-bifurcations type, limit cycles and center characterization} \label{se:preliminaries} In this section we present several technical results that are required for the statements and proofs developed in the next sections. More concretely, in Section~\ref{se:weakfocus} we outline the main features and differences of Hopf-type bifurcations in smooth and piecewise smooth settings, particularly regarding the bifurcation of limit cycles of small amplitude. How this bifurcation applies near centers is analyzed in Section~\ref{se:bifurcation}. A special case is the analysis of the return map near a monodromic point of $CC$-equilibrium type, which is addressed in Section~\ref{se:lyapunovconstants}. To assist the computations, we introduce a suitable small parameter and study the Taylor series of the return map near a center with respect to this parameter; this is done in Section~\ref{se:Lyapunovconstantsseries}. Since the classical center conditions are insufficient for our purposes in piecewise families, a new criterion is introduced in Section~\ref{se:sufficientcenter}. We finish with the analysis of the existence of a limit cycle from a monodromic equilibrium when a sliding segment appears and changes the stability. Previous results on unfoldings of this bifurcation do not apply here, as they fall outside our family. We require a new adapted result that guarantees the unfolding remains within our Kolmogorov class. This is developed in Section~\ref{se:pseudohopf}.

\begin{figure}[h]
	\begin{center}
		\begin{overpic}{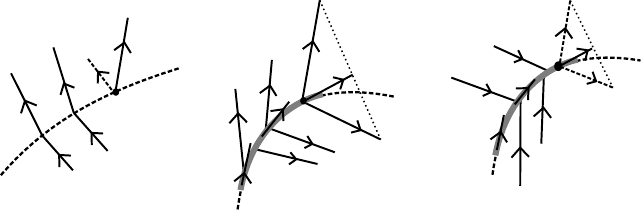}
			\put(29,20){$\Sigma$}
			\put(62.5,16){$\Sigma$}
			\put(101,21){$\Sigma$}
			\put(0,10){$\Sigma_1$}
			\put(32,10){$\Sigma_1$}
			\put(72,10){$\Sigma_1$}
			\put(2,2){$\Sigma_2$}
			\put(40,0){$\Sigma_2$}
			\put(83,2){$\Sigma_2$}
			\put(16,14){$p$}
			\put(46,14){$p$}
			\put(86,18){$p$}
			\put(9,30){$Z_1(p)=Z_i(p)$}
			\put(38,30){$Z_1(p)$}
			\put(94,14){$Z_1(p)$}
			\put(9,24){$Z_2(p)$}
			\put(57,7){$Z_2(p)$}
			\put(77,30){$Z_2(p)$}
			\put(56,21){$Z_i(p)$}
			\put(93,24){$Z_i(p)$}
		\end{overpic}
	\end{center}
	\caption{Definition of the vector field on $\Sigma$ following Filippov's convention in the crossing, escaping, and sliding regions.5}\label{fi:filipov}
\end{figure}
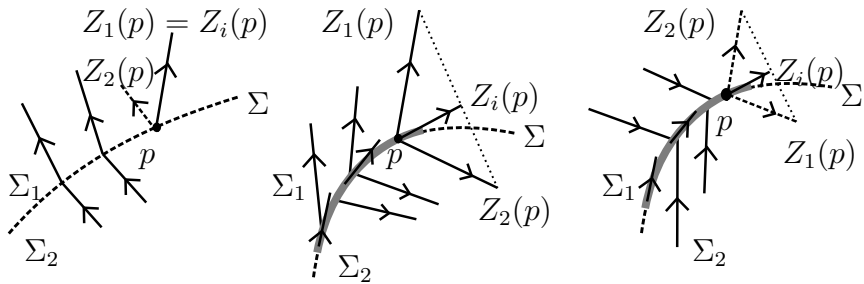

Let be a piecewise planar differential systems $Z=(Z_1,Z_2)$ and denote by $\Sigma$ its separation line. Points on the separation line where both $Z_1, Z_2$ simultaneously, point outward or inward from $\Sigma$ define the \emph{escaping} ($\Sigma^e$) and \emph{sliding region} ($\Sigma^s$), respectively. The interior of the complement of $\Sigma^e \cup \Sigma^s$ on $\Sigma$ defines the \emph{crossing region} ($\Sigma^c$). The boundaries of these regions are constituted by tangential points of $Z_1, Z_2$ with $\Sigma$. Even though in $\Sigma^e$ and $\Sigma^s$ the vector field $Z$ is multi-valuated, it is possible to define $Z$ by using the Filippov's convention, see Figure~\ref{fi:filipov}. Here we are interested only in the bifurcation analysis of crossing limit cycles and we do not present more details about the remaining concepts, for more details we refer the reader to \cite{Fil1988}.

\subsection{Return map near a monodromic point and weak focus order}\label{se:weakfocus}

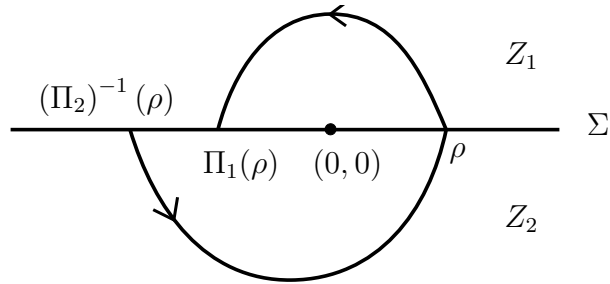
\begin{figure}[h]
	\begin{overpic}{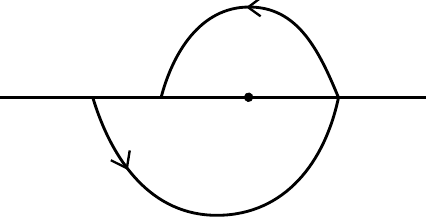}
		\put(55,20){$(0,0)$}
		\put(80,23){$\rho$}
		\put(35,20){$\Pi_1(\rho)$}
		\put(5,32){$\left(\Pi_2\right)^{-1}(\rho)$}
		\put(90,40){$Z_1$}
		\put(90,10){$Z_2$}
		\put(105,27){$\Sigma$}
	\end{overpic}
	\caption{The positive and negative half-return maps $\Pi_1$ and $\left(\Pi_2\right)^{-1},$ respectively.  }\label{fi:retunmaps}
\end{figure}

In smooth piecewise systems defined in two zones separated by a straight line, there are three distinct monodromic pseudo-equilibrium points: equilibrium-equilibrium, equilibrium-fold, and fold-fold. As in the smooth case, given a transversal section, the return map is defined in a neighbourhood of the pseudo-equilibrium point. For convenience, we use the separation straight line as the transversal section. Here, we will consider only the bifurcation of limit cycles from the first type. Let $Z = (Z_1, Z_2) \in \mathcal{K}_2$. Then the points $p \in \Sigma$ are simultaneously equilibria for both systems $Z_1$ and $Z_2$, but with a well-defined return map in a small neighbourhood, i.e., a monodromic behaviour. These are particular cases of the general concept of pseudo-equilibrium; see again \cite{Fil1988}. Moreover, we focus on the equilibrium-equilibrium points such that the Jacobian matrix of $Z_i$, $i=1,2$,  at each equilibrium has complex conjugated eigenvalues with nonzero imaginary part. Consequently, the return map can be obtained by the composition of the two half-return maps, denoted by $\Pi_1(\rho)$ and $\Pi_2 (\rho)$, associated to each component of $Z$. When $Z_1, Z_2$ are analytic then the return map also is. By simplicity, instead of using the composition of both maps, we will compute the difference map that is equivalent for our purposes. It is given by
\begin{equation}\label{eq:11}
	\Delta(\rho)=\left(\Pi_2\right)^{-1}(\rho)-\Pi_1(\rho)
\end{equation}
where, the function $\left(\Pi_2\right)^{-1}(\rho)$ is the inverse of the
negative half-return map $\Pi_2(\rho),$ as it is illustrated in
Figure~\ref{fi:retunmaps}. For  more details see  \cite{CollPhoGassu1999, GasTor2003}.

By analyticity, the return map, more precisely, each half-return map, can be expressed as a Taylor series near $\rho = 0$. As in the smooth scenario, the first nonvanishing coefficient of the Taylor series of $\Delta$ is called the $k$th-order Lyapunov quantity. It is denoted by $W_k$ and is defined only modulo $\{W_1 = \cdots = W_{k-1} = 0\}$. 

In the study of smooth planar vector fields, the first nonvanishing coefficient of the return map always has an odd subscript: $L_{\ell}=W_{2\ell+1} \ne 0$. It can be seen that, generically, the number of limit cycles that bifurcate from a monodromic point is at most $\ell$. In this case, we say that the weak focus order is also $\ell$, because the order is related to the number of limit cycles that a versal unfolding has using a generic analytic perturbation. The usual approach to study this bifurcation is to restrict the analysis to perturbations without constant or linear terms. In this case, $W_1 = W_2 = 0$, and the coefficients $W_k$ are polynomials in the perturbation parameters. Generically, $L_1=W_3$ is nonvanishing and determines the stability of the origin. A limit cycle of small amplitude may appear when the stability changes due to the inclusion of the trace parameter. This is the well-known codimension-$1$ Hopf bifurcation. More limit cycles can appear when this first coefficient vanishes, in what is called a degenerate Hopf bifurcation, which has codimension higher than $1$. In this case, the Bautin ideal $\mathcal{B} = \langle W_3, \ldots, W_k, \ldots \rangle$ is introduced. This ideal has a finite number of generators when the perturbation class is polynomial, and it is generated only by the elements with odd subscripts. See, for example, \cite{CimGasMan2020, LiuLiHua2014, RomaSha2009}. This property implies that the number elements in $\mathcal{B}$ decreases by half, and when the bifurcation is versal $\ell$ limit cycles appear from a weak focus of order $\ell$, in this case there exist perturbation parameters such that the elements in the Bautin ideal alternate in sign, or, alternatively, are ``independent.'' Although this bifurcation is discussed in classical references such as \cite{AndLeoGorMai1973, Rou1998}, a more detailed analysis will be presented later.

The same phenomenon is observed in the smooth piecewise scenario when the equilibrium point is of fold-fold type (i.e., when both vector fields have coincident invisible tangency points on $\Sigma$). In this case, however, the relevant quantities have even subscripts, with the first nonvanishing one determining the stability behavior. Hence, when $L_\ell = W_{2\ell} \ne 0$, we say that the origin has a weak focus of order $\ell$. Generically, as in the classical high-codimension Hopf bifurcation, $\ell$ crossing-type limit cycles bifurcate from the origin under a generic perturbation. In the literature, this equivalent bifurcation is known as the pseudo-Hopf bifurcation. It was first discovered by Filippov in \cite{Fil1988}, and the term was coined in \cite{KuzRinGra2003}. See \cite{EstebanPonce2022, NovaesLeandro2021} for more details on the higher-codimension case.

The special symmetry present in the above cases is not satisfied in a generic piecewise perturbation near a monodromic point. A monodromic equilibrium $p \in \Sigma$ is said to be a weak focus of order $\ell$ when $W_j = 0$ for $1 \leq j \leq \ell - 1$ and $W_\ell \ne 0$. Its stability, as usual, is determined by the sign of the first nonzero Lyapunov quantity $W_\ell$, and under a generic perturbation, $\ell$ small-amplitude limit cycles bifurcate from the equilibrium point, which can be taken at the origin. This bifurcation is also known as the degenerate or higher-codimension pseudo-Hopf bifurcation, and the existence of a generic versal unfolding is established in \cite{GouvTorre2021l}. To apply it within our Kolmogorov family, the analysis of the bifurcation of the limit cycle associated with the appearance of a sliding segment requires a refinement; see Section~\ref{se:pseudohopf}.

\subsection{The bifurcation of limit cycles from a center family in piecewise systems}\label{se:bifurcation}

Christopher \cite{Chr2005} explains how the $k$-th order Taylor approximation of the Lyapunov quantities can be used to obtain lower bounds on the number of limit cycles near center-type equilibrium points in smooth systems. These bounds are related to the intersection of algebraic varieties and the transversality properties of Taylor expansions, either at the linear level or at higher orders. Similar results can be found in \cite{ChiJac1989, ChiJac1991, Han1999}.

\begin{pro}[\cite{Chr2005}]\label{thm:chrorder1}
Suppose that $s$ is a point on the center variety and that the first $k$ Lyapunov quantities, $W_{1},\ldots,W_{k},$ have independent linear parts (with respect to the expansion of $W_{j}$ about $s$), then $s$ lies on a component of the center variety of codimension at least $k$ and there are bifurcations which produce $k$ limit cycles locally from the center-type equilibrium point corresponding to the parameter value $s$. If, furthermore, we know that $s$ lies on a component of the center variety of codimension $k$, then $s$ is a smooth point of the variety, and the cyclicity of the center for the parameter value $s$ is exactly $k$. In the latter case, $k$ is also the cyclicity of a generic point on this component of the center variety.
\end{pro} 

The subsequent proposition, also due to Christopher \cite{Chr2005}, is an extension of the above result which shows that sometimes we can obtain more limit cycles using high-order Taylor developments of the Lyapunov quantities. 

\begin{pro}[\cite{Chr2005}]
Suppose that we are in a point $s$ where Proposition~\ref{thm:chrorder1} applies. After a change of variables if necessary, we can assume that $W_1=W_2=\cdots=W_k=0$ and the next Lyapunov quantities $W_{j}=h_j(\lambda)+O_{m+1}(\lambda),$ for $j=k+1,\ldots,k+l,$ where $h_j$ are homogeneous polynomials of degree $m\ge2$ and $\lambda=(\lambda_{k+1},\ldots,\lambda_{k+l})$. If there exists a line $\mathcal{L},$ in the parameter space, such that $h_{j}(\mathcal{L}) = 0,$ $j=k+1,\ldots,k+l-1,$ the hypersurfaces $h_{i} = 0$ intersect transversally along $\mathcal{L}$ for $j=k+1,\dots,k+l-1,$ and $h_{k+l}(\mathcal{L})\ne 0$, then there are perturbations of the center which can produce $k+l$ limit cycles.
\end{pro}

The same method can be applied to perturbations of a family of piecewise polynomial systems with a monodromic pseudo-equilibrium point of center type. Proposition~\ref{prop:theop} is an extended version of Proposition~\ref{thm:chrorder1}, originally due to Gouveia and Torregrosa in \cite{TorreGouv2021h}. Before stating it, we introduce some notation.

Following the notation introduced in this paper, we consider two regions separated by a straight line passing through the origin, each denoted by $i = 1, 2$. We study the piecewise family of vector fields $Z = (Z_{1,\mu}, Z_{2,\mu})$, defined by the solutions of the differential equations $(\dot{x}, \dot{y}) = (P^i_c(x, y, \mu), Q^i_c(x, y, \mu))$, where each system is polynomial of degree $n$, depending on a parameter $\mu \in \mathbb{R}^\ell$, and has a center equilibrium point at the origin. That is, we consider the perturbed polynomial system
\begin{equation}\label{eq:centrob}
	\begin{array}{lll}
		\dot{x} & = & X_{i,c}(x, y, \mu) + P_i(x, y, \lambda), \\
		\dot{y} & = & Y_{i,c}(x, y, \mu) + Q_i(x, y, \lambda),
	\end{array}
\end{equation}
where $P_i$ and $Q_i$ are polynomials of degree $n$. In particular,
\begin{equation*}
	P_i(x, y, \lambda) = \sum\limits_{k+l=0}^n a^i_{k,l} x^k y^l, \quad Q_i(x, y, \lambda) = \sum\limits_{k+l=0}^n b^i_{k,l} x^k y^l,
\end{equation*}
with $\lambda = (a^i_{00}, a^i_{10}, a^i_{01}, \ldots, b^i_{00}, b^i_{10}, b^i_{01}, \ldots) \in \mathbb{R}^{M}$, where $M = 2n^2 + 6n + 4$. We denote $X_c = X_{i,c}$ and $Y_c = Y_{i,c}$ when the centers are the same on both sides. The next proposition holds in this context.

\begin{pro} [\cite{TorreGouv2021h}] \label{prop:theop}	We denote by $W_{j}^{[1]}(\lambda,\mu)$ the first-order development, with respect to $\lambda\in\mathbb{R}^M,$ of the $j$-Lyapunov quantity of system~\eqref{eq:centrob}, for each fixed value of $\mu\in\mathbb{R}^\ell$. We assume that, after a change of variables in the parameter space if necessary, we can write
\begin{equation*}
	W_j=\begin{cases}
		&\lambda_j + O_2(\lambda), \text{ for } j=1,\ldots,k-1,\\
		&\sum\limits_{l=1}^{k-1} g_{j,l}(\mu) \lambda_l+f_{j-k}(\mu)\lambda_{k}+ O_2(\lambda), \text{ for } j=k,\ldots,k+l.
	\end{cases}
\end{equation*}
where with $O_2(\lambda)$ we denote all the monomials of degree higher or equal than $2$ in $\lambda$ with coefficients analytic functions in $\mu$. If there exists a point $\mu^*$ such that $f_0(\mu^*)=\cdots=f_{l-1}(\mu^*)=0,$ $f_{l}(\mu^*)\ne0,$ and the Jacobian matrix of $(f_{0},\ldots,f_{l-1})$ with respect to $\mu$ has rank $l$ at $\mu^*,$ then system~\eqref{eq:centrob} has $k+l$ hyperbolic limit cycles of small amplitude bifurcating from the origin. 
\end{pro}
Usually, as previously described, in higher-codimension Hopf-type bifurcations the strategy to obtain small-amplitude limit cycles begins by analyzing the perturbation without linear or constant terms. After identifying the resulting limit cycles, we obtain an additional one by including the linear terms, changing the sign of the trace while keeping the equilibrium at the origin, and finally another by including the constant terms, which affect stability through the sliding segment. All of these are taken into account in the result above.

As noted in \cite{Rou1998}, the main obstruction to obtaining an unfolding of limit cycles near a Hopf point lies in the difficulty of analyzing the independence of the Lyapunov quantities, or in ensuring that they alternate in sign. The results of this section offer a computational approach to reach such conclusions near a center point. Similar results can be established near a weak focus point, since the methodology is primarily based on the intersection analysis of polynomial varieties in the parameter space.

\subsection{The difference map in piecewise systems near a $CC$-type equilibrium point}\label{se:lyapunovconstants} 

Here we present how to obtain the coefficients $W_i$ of the difference map \eqref{eq:11} near a monodromic pseudo-equilibrium, $(x_0,y_0)$, of a piecewise polynomial system $Z=(Z_1, Z_2)$ of degree $n$, where $\Sigma$ is a straight line passing through $(x_0,y_0)$, using first integrals of $Z_1$ and $Z_2$.

We can assume, after an adequate affine change of variables if necessary, that $(x_0,y_0)=(0,0)$ and $\Sigma$ is the $x$-axis. For our purposes, we will assume that each $Z_i$ has a first integral $H_i:\mathbb{R}^2\rightarrow\mathbb{R}$, for $i=1,2$, which we will assume to be analytic and whose Taylor series are written as $H_i(x,y)=A_i x^{2}+B_i x y+C_i y^2+O_3(x,y)$. As we are analyzing the solutions near $y=0$, the monodromy condition guarantees that $A_i \ne 0$, and by dividing the first integrals by $A_i$, we can take $A_i=1$.

Let $(\rho,0)$ and $(\sigma,0)$ be two points defined by the intersection of a level curve of the first integral $H_1$ with the separation straight line $y=0$. Then, they satisfy the equation $H_1(\rho,0)-H_1(\sigma,0)=(\rho-\sigma)\mathcal{H}_1(\rho,\sigma)=0$ for an analytic function $\mathcal{H}_1(\rho,\sigma)$.

As the Taylor series of $\mathcal{H}_1(\rho,\sigma)$ near $(0,0)$ is written as $\rho+\sigma+O_2(\rho,\sigma)$, the Implicit Function Theorem allows us to write $\sigma$ as a function of $\rho$. In fact, this is the half-return map $\Pi_1$, defined previously. Moreover, we can obtain, in a recursive way, the series
\[
\Pi_1(\rho)=-\rho+\sum_{k=2}^{\infty} W_{1,k} \rho^k.
\]

With the same technique, we can obtain the value of $\sigma$ for the first integral $H_2$, which gives the inverse of the function $\Pi_2$, yielding the series
\[
(\Pi_2)^{-1}(\rho)=-\rho+\sum_{k=2}^{\infty} W_{2,k} \rho^k,
\]
from the function $\mathcal{H}_2(\rho,\sigma)$. Hence, the difference map \eqref{eq:11} is
\begin{equation}\label{eq:difmap}
	\Delta(\rho)=(\Pi_2)^{-1}(\rho)-\Pi_1(\rho)=-\sum\limits_{k=2}^\infty(W_{2,k}
	-W_{1,k})\rho^k=-\sum\limits_{k=2}^\infty W_k\rho^k.
\end{equation}

In the following, we present the expressions of $W_1,\ldots,W_8$. These coefficients can be determined from the Taylor series
\[
H_i(\rho,0)=\rho^2+\sum\limits_{k=3}^\infty h_{i,k}\rho^k.
\]

Then, using the mechanism described above, we obtain the coefficients of the half-return maps:
\begin{equation*}
	\begin{aligned}
		W_{i,2}&= -h_{i,3},\\
		W_{i,3}&= -h_{i,3}^2,\\
		W_{i,4}&= -2 h_{i,3}^3+2 h_{i,3} h_{i,4}-h_{i,5},\\
		W_{i,5}&= -4 h_{i,3}^4+6 h_{i,3}^2 h_{i,4}-3 h_{i,3} h_{i,5},\\
		W_{i,6}&= -9 h_{i,3}^5+19 h_{i,3}^3 h_{i,4}-11 h_{i,3}^2 h_{i,5}-4 h_{i,3} h_{i,4}^2+3 h_{i,3} h_{i,6}+2 h_{i,4} h_{i,5}-h_{i,7},\\
		W_{i,7}&= -21 h_{i,3}^6+56 h_{i,3}^4 h_{i,4}-34 h_{i,3}^3 h_{i,5}-24 h_{i,3}^2 h_{i,4}^2+12 h_{i,3}^2 h_{i,6}+16 h_{i,3} h_{i,4} h_{i,5}\\
		&\quad -4 h_{i,3} h_{i,7}-2 h_{i,5}^2,\\
		W_{i,8}&=
		-51 h_{i,3}^7+165 h_{i,3}^5 h_{i,4}-104 h_{i,3}^4 h_{i,5}-112 h_{i,3}^3 h_{i,4}^2+43 h_{i,3}^3 h_{i,6}\\
		&\quad +93 h_{i,3}^2 h_{i,4} h_{i,5}+8 h_{i,3} h_{i,4}^3-18 h_{i,3}^2 h_{i,7}-12 h_{i,3} h_{i,4} h_{i,6}-17 h_{i,3} h_{i,5}^2\\
		&\quad -4 h_{i,4}^2 h_{i,5}+4 h_{i,3} h_{i,8}+2 h_{i,4} h_{i,7}+3 h_{i,5} h_{i,6}-h_{i,9}.
	\end{aligned}
\end{equation*}

In this case, for a general problem, we obtain the first nonvanishing coefficients of the difference map \eqref{eq:11}. We notice that, by construction, $W_1=0$ and, after straightforward simplifications, for $k=2,\ldots,8$ and under the conditions $\{W_1=W_2=\cdots=W_{k-1}=0\}$, we get $W_3=W_5=W_7=0$ and
\begin{equation*}
	\begin{aligned}
		W_2&=-h_{1,3}+h_{2,3},\\
		W_4&=2 h_{1,3} h_{1,4}-2 h_{1,3} h_{2,4}-h_{1,5}+h_{2,5},\\
		W_6&=-3 h_{1,3}^3 h_{1,4}+3 h_{1,3}^3 h_{2,4}-4 h_{1,3} h_{1,4}^2+4 h_{1,3} h_{1,4} h_{2,4}+3 h_{1,3} h_{1,6}\\
		&\quad -3 h_{1,3} h_{2,6}+2 h_{1,4} h_{1,5}-2 h_{1,5} h_{2,4}-h_{1,7}+h_{2,7},\\
		W_8&=11 h_{1,3}^5 h_{1,4}-11 h_{1,3}^5 h_{2,4}+28 h_{1,3}^3 h_{1,4}^2-28 h_{1,3}^3 h_{1,4} h_{2,4}-11 h_{1,3}^3 h_{1,6}+11 h_{1,3}^3 h_{2,6}\\
		&\quad -11 h_{1,3}^2 h_{1,4} h_{1,5}+11 h_{1,3}^2 h_{1,5} h_{2,4}+8 h_{1,3} h_{1,4}^3-8 h_{1,3} h_{1,4}^2 h_{2,4}-12 h_{1,3} h_{1,4} h_{1,6}\\
		&\quad +6 h_{1,3} h_{1,4} h_{2,6}+6 h_{1,3} h_{1,6} h_{2,4}-4 h_{1,4}^2 h_{1,5}+4 h_{1,4} h_{1,5} h_{2,4}+4 h_{1,3} h_{1,8}\\
		&\quad -4 h_{1,3} h_{2,8}+2 h_{1,4} h_{1,7}+3 h_{1,5} h_{1,6}-3 h_{1,5} h_{2,6}-2 h_{1,7} h_{2,4}-h_{1,9}+h_{2,9}.
	\end{aligned}
\end{equation*}

We remark that in the above simplification one can see the same phenomenon of vanishing half-terms in the Taylor developments as in the fold-fold equilibrium. Consequently, after the center problem is solved, the number of limit cycles that will bifurcate from the origin will depend on the perturbation. If the perturbation vector fields are centers on both sides, only even terms are useful and the bifurcation mechanism of the fold-fold equilibria applies. Note that the above expressions of $W_{2i}$ are independent because of the linear terms $-h_{1,2i+1}+h_{2,2i+1}$, but, as we will see, this is not the case when we fix the perturbation family. Hence, the unfolding analysis is more intricate. If the perturbation is generic, as we will see in a later section, then the above expressions are not enough to study the difference map, because the odd terms also appear and, consequently, more limit cycles can arise.

A natural consequence of this analysis is that when we consider piecewise vector fields having centers on both sides, the expected number of limit cycles is lower, in fact, less than half, than in the general case. This fact can be seen in the piecewise linear class, where the total number found so far is three, as we have mentioned in the introduction, but no limit cycles of crossing type appear when both sides have a center. A very simple proof of this fact follows from the algorithm described in this section. The first integrals of a linear system having a center are polynomials of degree two in $x,y$. In this case, it is not restrictive to assume that the separation line is $y=0$. Then the return map follows from $\mathcal{H}_1(\rho,\sigma)=a \rho + a\sigma + b$ for some real constants $a,b$. Then, as the half-return map is written as $\sigma=-\rho-b/a$, the difference map is constant; consequently, it has no zeros.

\subsection{Taylor series of the return map with respect to a privileged small parameter}\label{se:Lyapunovconstantsseries} 
Perturbation theory usually requires, more or less explicitly, the knowledge of a solution. Here, we briefly recall the algorithm described in more detail in \cite{BraCruTor2024} for computing the Taylor series with respect to $\varepsilon$ of a solution which, for $\varepsilon = 0$, corresponds to a center. This algorithm is inspired by the one introduced in \cite{GasTor2003}.
Assume that $Z=(Z_{1,\mu}, Z_{2,\mu})$ is a piecewise polynomial system having the origin as a monodromic pseudo-equilibrium in which there exists a neighbourhood where all the trajectories are periodic for any value of the parameter $\mu$, for simplicity we say that this is a monodromic pseudo-equilibrium of center type. Denotes by $Z_{i,\varepsilon}$ a piecewise perturbation as in \eqref{eq:centrob} of degree $n$ of $Z_{i,\mu}$, then
\begin{equation}\label{eq:a8}
	Z_{i,\varepsilon}:(\dot{x},\dot{y})=Z_{i,\mu}(x,y)+\varepsilon Z_{i}(x,y),  \textrm{ if } (x,y)\in \Sigma_i,
\end{equation}
with  $\varepsilon > 0,$ $Z_{i,\mu}=(X_{i,c}(x,y,\mu),Y_{i,c}(x,y,\mu))$  and $Z_{i}(x,y) = \left(X_{i}(x,y,\lambda),Y_{i}(x,y,\lambda)\right)$, $X_i$, $Y_i$ are polynomial vector fields of degree $n$ without constant terms, for $i=1,2$. Here, we also denote $X_c=X_{i,c},$ and $Y_c=Y_{i,c}$ when the centers are the same in both side.

Applying polar coordinates $(x,y)=(r\cos\theta,r\sin\theta)$, in the perturbed center \eqref{eq:a8}, and making $\theta$ as the new independent variable, we get the solution of the differential equation associated to \eqref{eq:a8}, 
$$\dot{r}=\dfrac{dr}{d\theta}=F(\theta, r)=F_0(\theta, r) + \sum_{j=1}^m\varepsilon^j F_j(\theta, r) + \mathcal O(\varepsilon^{m+1}),$$ in which $F_j=(F_{1,j}, F_{2, j}) $, $j=0,\ldots, m$. 
Such a solution can be written as $\varphi(\theta,r)=(\varphi_{1,N}(\theta,r), \varphi_{2,N}(\theta,r))$ where $\varphi_{i,N}(\theta,r)=\sum_{k=0}^N\varepsilon^k\varphi_{i,k}(\theta,r),$  
for $i=1,2,$ and $N$ a natural number. 
Thus, we define the \emph{$k$-difference function} as  
\begin{equation*}
	\delta_k(r) = \varphi_{1,k}(\pi,r) - \varphi_{2,k}(-\pi,r), 
\end{equation*}
where $\varphi_{i,0}$, $i=1,2,$ is the solution of the initial value problem 
$$z'(s) = F_{i,0}(s,z(s)),\quad z(0) = r,$$ and $\varphi_{i,k}$'s, for $k\geq1$, are given recursively adapting \cite{LliNovTei2014}. Then, from this we can defined the $N$-jet of the \emph{difference function} as 
\begin{equation}\label{delta1}
	\Delta_N(r,\varepsilon)=\sum_{k=1}^N \varepsilon^k \delta_k(r).
\end{equation}
Observe that, as the origin of the non-perturbed vector field is of center type,
$\varphi_{1,0}(\pi,r) - \varphi_{2,0}(-\pi,r) = 0$, so $\delta_0$ does not appear. 
Clearly, each simple zero of $r\mapsto \Delta(r, \varepsilon)$ provides a hyperbolic limit cycle of \eqref{eq:a8}. To compute the $k$-th order Lyapunov quantities, we propose another approach: instead of looking for explicit formulas for $\delta_l(r)$, we consider their Taylor expansions in $r$. This leads to an algorithm that only requires integrating trigonometric functions. Moreover, we also propose a type of \emph{blowing up} technique  in such a way that the $N$-jet of the difference functions \eqref{delta1}, is written as follows
\[
\Delta_N(r)=\sum_{j=1}^\infty \mathcal{W}^{[N]}_{j} r^j,
\] 
where $\mathcal{W}^{[N]}_{j}=\sum_{k=1}^N W^{[k]}_{j},$ with $W^{[k]}_{j}$ is the homogeneous polynomial of degree  $k$ of the $j$-th Lyapunov quantity in the coefficients $\lambda.$ Since each coefficient $W_j$ has an analytic dependence with respect to the perturbation parameters. So, we can conclude that the $j$-th Lyapunov quantity, with a slight abuse of notation, as $W_j=\mathcal{W}^{[\infty]}_{j}=\sum_{k=1}^{\infty} W^{[k]}_{j},$ and the complete difference functions as   
\[
\Delta(r)=\sum_{j=1}^\infty W_j r^j.
\] 
Henceforth, we will write either $W_j^{[1]}$ or $\mathcal{W}_j^{[1]}$ to refer to the same object.
 
\subsection{Sufficient conditions to be a center in piecewise systems}\label{se:sufficientcenter}

For smooth systems, if the equilibrium point is of monodromic type and there is a first integral, then the equilibrium point is a center. However, for non-smooth systems, we must prove that the positive and negative half-return maps satisfy $\left(\Pi_2\right)^{-1}(\rho)-\Pi_1(\rho)=0$ for any $\rho>0$, see \eqref{eq:11} and Figure~\ref{fi:retunmaps}. Different situations about how to check this condition can be seen, for example, in \cite{daCruzTor2022}. In particular, getting the first integrals in each region is not enough to have a local first integral well defined and continuous in an open set. So we will say that a piecewise system $Z=(Z_1, Z_2)$ has a \textit{$\Sigma$-first integral }if $Z_1$ and $Z_2$ have first integrals in each region, $\Sigma_1$ and $\Sigma_2$, respectively. The usual definition of first integral, that is, a non-constant function that is constant along the solutions implies the continuity condition. Consequently, if we have a continuous piecewise first integral around a monodromic pseudo-equilibrium we will have a center. This situation will be enough for this paper.

Here, we take the $x$-axis as the separation line to simplify the reading in the following result. The result can be easily generalized, considering other separation lines. For example, another straight line, when necessary, is enough to apply a rotation. 

\begin{pro}\label{teocentros} Let $Z=(Z_1, Z_2)$ be a planar piecewise differential system with a monodromic pseudo-equilibrium at the origin and $\Sigma=\{y=0\}$. If there exist first integrals $H_i$ of $Z_i$ satisfying $H_1(x,0)=H_2(x,0),$  where $i=1,2.$ Then, $Z$ has a center at the origin.  
\end{pro}   	
	
\subsection{The pseudo-Hopf bifurcation form }\label{se:pseudohopf}  
\begin{figure}[ht!]
	\includegraphics[height=1.5cm]{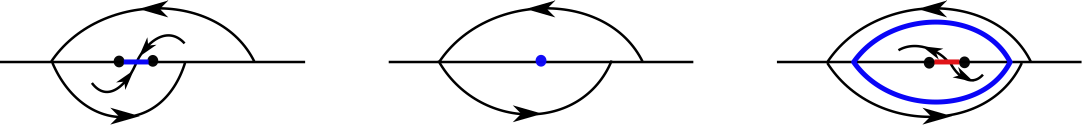}
	\caption{The Hopf-type bifurcation}\label{fi:hopftype1}
\end{figure}

\begin{figure}[ht!]
	\includegraphics[height=8cm]{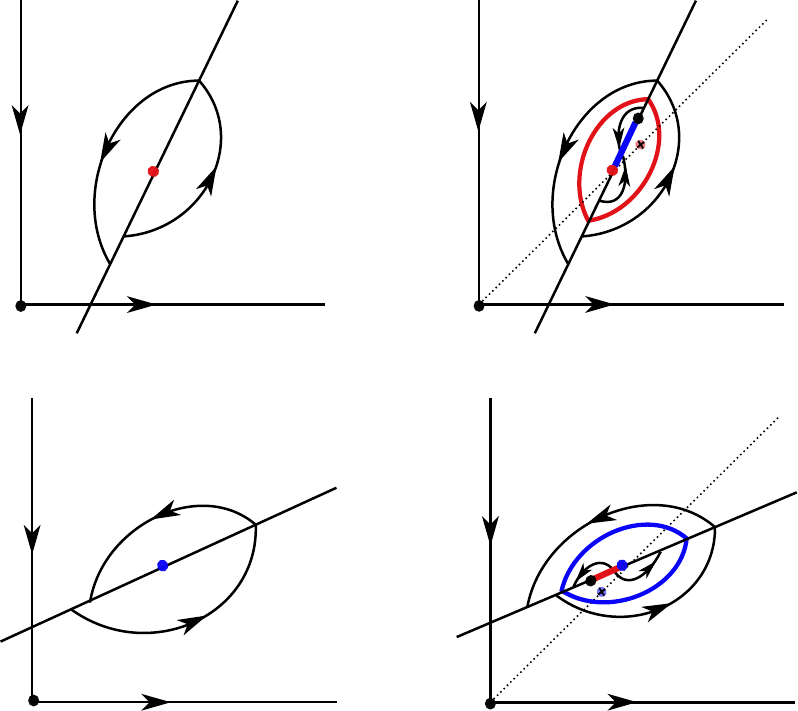}
	\caption{Some Hopf-type bifurcations. The invisible equilibrium is drawn in the doted straight line $y-x=0$}\label{fi:hopftype}
\end{figure}

As previously mentioned, in piecewise differential systems, a limit cycle can arise from a monodromic pseudo-equilibrium, since the stability of the point can be altered by adding a sliding or escaping segment. This phenomenon for fold-fold equilibria was termed a pseudo-Hopf bifurcation in \cite{KuzRinGra2003}, and it had already been proven in \cite{Fil1988}. See Figure~\ref{fi:hopftype1}. A collection of similar Hopf-type bifurcations can be found in \cite{Sim2022}. The next Theorem~\ref{hopf-type} ensures this fact in a general context, as in \cite{BraCruTor2024, FrePonTor2014}. However, as we need a perturbation within the Kolmogorov class, a new result needs to be stated, because the perturbations must preserve the invariant lines $x = y = 0$. We will see that the stability of the monodromic pseudo-equilibrium with respect to the relative position of the separation straight line is essential to guarantee this bifurcation of a limit cycle. The way to overcome this obstacle in order to obtain a Hopf-type bifurcation can be seen in the proof of Proposition~\ref{hopf_bif}. In all cases, the perturbed system has no visible equilibria outside the separation line, and the limit cycles inherit the stability from the equilibrium. An illustration is shown in Figure~\ref{fi:hopftype}.

\begin{tho}\label{hopf-type} Let $Z=(Z_1, Z_2)$ be a planar piecewise differential system and $(x_0,y_0)\in\Sigma$ a monodromic pseudo-equilibrium point, stable (resp. unstable). Then, if an escaping segment (resp. sliding segment) appears in the neighbourhood of $(x_0,y_0)\in\Sigma,$ we obtain a stable (resp. unstable) limit cycle, and this phenomenon is called pseudo-Hopf or Hopf-type bifurcation.  
\end{tho}

\begin{proof}
The proof is a direct consequence of the Poincaré--Bendixon Theorem for piecewise differential systems, see \cite{BuzCarEuz2018}. Directly, having a monodromic pseudo-equilibrium point stable, and doing an  repulsing escaping segment over $\Sigma$ appears, changing the stability of the neighbourhood of the point $p\in\Sigma$ and an stable limit cycle of small amplitude bifurcates. We observe that choosing an suitable perturbation of a monodromic pseudo-equilibrium an escaping (or sliding) segment appears.
\end{proof}

\begin{pro}\label{hopf_bif} Let $Z = (Z_1, Z_2)$ be the piecewise Kolmogorov system defined by $Z_i := (\dot{x}, \dot{y}) = (x f_i(x,y), y g_i(x,y))$ for $i = 1,2$ with the separation straight line $\alpha (x - 1) - \beta (y - 1)=0$ being $(\alpha, \beta) \neq (0,0),$ $\beta>0$, and $\alpha+\beta\ne0$. Let $(1,1)$ be a monodromic pseudo-equilibrium such that the vector field rotates counter-clockwise. Then, if $|\alpha|\ge \beta$ and the equilibrium is unstable (resp. $|\alpha|\le\beta$ and stable) there exists a perturbation inside Kolmogorov class such that a unstable (resp. stable) limit cycle bifurcates from it.
\end{pro}	

\begin{proof} For each $i=1,2$ we consider the perturbed system $Z_{i,\varepsilon}$ of $Z_i$ obtained under the homothetic change $(x,y) \rightarrow ((1+\varepsilon)x, (1+\varepsilon)y)$. Fixing $\alpha,\beta,$ the proof follows by choosing an appropriate sign for $\varepsilon$ and perturbing only one vector field. When $\alpha/\beta \geq 1,$ as $(1,1)$ is unstable, the vector field $(Z_{1,\varepsilon}, Z_2)$ exhibits a Hopf-type bifurcation, creating an unstable limit cycle of small amplitude around $(1,1)$ for small $\varepsilon > 0$. Similarly, when $0 \leq \alpha/\beta \leq 1,$ the equilibrium is stable, and the limit cycle, which is also stable, bifurcates when $\varepsilon < 0$. The proof for the other conditions follows similarly, taking $(Z_1,Z_{2,\varepsilon})$ with $\varepsilon>0$.
\end{proof}
	
We note that the zone being changed and the sign of $\varepsilon$ depend on the homothetic transformation we have chosen with respect to the separation line.

\section{Caracterization of some classes of centers in $\mathcal K_2$} \label{centers}
In this section, we will work on the characterization of the center problem for some classes of $\mathcal{K}_2$. As mentioned before, we restrict our study to the case where the pseudo-equilibrium is a $CC$-equilibrium point and some specific separation straight lines. Characterizing centers for every separation straight line is a very intricate task. To show the difficulty of this problem, in Proposition~\ref{center2}, we solve the center problem by fixing one parameter of the system and selecting a specific straight line of separation. As we have commented above, having two separate centers whose equilibria coincide is not enough to guarantee a coupling center of the original system. Due to the computational difficulties, our aim is not to provide a complete classification of the center problem; we can think of this as a first attempt at a generalization of the center problem for the classical Lotka--Volterra system. As we have also explained before, a system of type
\begin{equation}\label{eq:LV}
	\begin{cases}
		\dot{x}=x (a+b x+c y),\\
		\dot{y}=y (d+e x+f y),
	\end{cases}
\end{equation}
having an equilibrium point in the first quadrant of monodromic type with zero trace of the Jacobian matrix is a center. We recall that the classical Lotka--Volterra system is a particular case of \eqref{eq:LV} taking $b=f=0$, for which the first integral is implicitly obtained from a direct integration using the separation variables method. We will recover this system in Section~\ref{bifurcation}, where it will be perturbed in the piecewise class $\mathcal{K}_{2}$. Obviously, the center property will not depend on the separation line in this case. Consequently, it will appear in the families described in Theorem~\ref{c3} and Proposition~\ref{center2}.

\begin{proof}[Proof of Theorem~\ref{c3}]
System \eqref{eq:2}, when $b_i f_i - c_i e_i \ne 0$ for $i=1,2$, has an equilibrium point outside the coordinate axes, $x_0 y_0 \ne 0$, located at
\[
(x_0, y_0) = \left( \frac{c_i d_i - a_i f_i}{b_i f_i - c_i e_i}, \frac{a_i e_i - b_i d_i}{b_i f_i - c_i e_i} \right), \text{ for } i=1,2.
\]
We can assume, interchanging $x$ with $y$ if necessary, that $e_i \ne 0$. Moreover, through straightforward computations, we obtain the necessary conditions 
\[
D_i = 2b_i e_i x_0^2 + 4b_i f_i x_0 y_0 + 2c_i f_i y_0^2 + (a_i e_i + 2b_i d_i)x_0 + (2a_i f_i + c_i d_i)y_0 + a_i d_i > 0,
\]
with
\[
a_i = \frac{(b_i x_0)^2 - b_i e_i x_0^2 + D_i^2}{e_i x_0},\quad
c_i = -\frac{(b_i x_0)^2 + D_i^2}{e_i x_0 y_0},\quad
d_i = (b_i - e_i) x_0,\quad
f_i = -\frac{b_i x_0}{y_0},
\]
for $i=1,2$, under which $(x_0, y_0)$ is a $CC$-equilibrium point, that is, a point whose Jacobian matrices have zero trace and positive determinant. Although we could write the center conditions, detailed in the statement, depending only on the parameters of $Z$, due to their size, we have chosen more compact expressions. Moreover, to simplify computations, we can also assume, rescaling and reparametrizing time if necessary, that the $CC$-equilibrium point is located at $(1,1)$ and $D_i = 1$ for $i=1,2$.
	
Under the above hypotheses, each system $Z_i$, for $i=1,2$, has a first integral of the form
\begin{equation}\label{ip}
H_i(x,y) = x^{b_i(b_i - e_i)} y^{-\frac{b_i}{e_i}(b_i^2 - b_i e_i + 1)} \Lambda_i(x,y)
\end{equation}
where
\begin{equation*}
\Lambda_i(x,y) = (b_i^2 e_i - b_i e_i^2)x + (b_i^2 e_i - b_i^3 - b_i)y + b_i^3 - 2b_i^2 e_i + b_i e_i^2 + b_i - e_i.
\end{equation*}
Then, after a translation, the equilibrium point is located at the origin and the algorithm to compute the $n$-th Lyapunov constant of $Z$, detailed in Section~\ref{se:lyapunovconstants}, applies after adding the monodromic condition. We note that, in some cases, an extra rotation is necessary.
	
The strategy for proving all cases is the same. Straightforward computations provide the Lyapunov quantities. We note that the vanishing trace condition guarantees that $W_1 = 0$. Then, after removing a nonvanishing factor depending on the parameters, we have that $W_2$ only vanishes under each condition of the statement. Moreover, $W_3 = W_4 = \cdots = W_8 = 0$ when $W_2 = 0$. So, we have only one necessary condition for having a center. The sufficient condition follows from Proposition~\ref{teocentros}, checking the continuity of expressions \eqref{ip} defining each first integral. We remark that, after performing the inverse of the rescaling transformations, we recover the expression of the center condition given in the statement.
\end{proof}

Our second approach involves selecting a straight line different from those in Theorem~\ref{c3}. This adjustment rendered the problem even more challenging,  the first Lyapunov constant is quite bigger, leading to a further increase in difficulty, as evidenced in our subsequent result. Therefore, we opted to restrict the center-focus problem even more, fixing one of the parameters.

\begin{pro}\label{center2}  Let  $Z\in \mathcal K_2$ with a $CC$-equilibrium point on $\Sigma$ and $b_1=1$. Assuming that the separation straight line is $4(x-1)-3(y-1)=0,$ $Z$ has a center at the $CC$-equilibrium point if, and only if, one of the next condition hold:
	\begin{itemize}
		\item[$C_1:$] $b_2-1=e_1-e_2=0;$
		\item[$C_2:$] $b_2+1=e_1+e_2=0;$	
		\item[$C_3:$] $3e_1-8=3e_2-4b_2=0;$
		\item[$C_4:$] $3e_1-8=3e_2b_2-4 \left(b_2^2+1\right)=0;$ 
		\item[$C_5:$] $3e_1-8=(8b_2-7e_2)^2-49(b_2^2-64)=0;$ 
		\item[$C_6:$] $3e_1-4=3e_2-4b_2=0;$
		\item[$C_7:$] $3e_1-4=3e_2b_2-4(b_2^2+1)=0;$
		\item[$C_8:$] $3e_1-4=(8b_2-7e_2)^2-49(b_2^2-64)=0.$
	\end{itemize} 
\end{pro}

\begin{proof}
Doing as in the proof of the previous result, we compute the first four (nonvanishing) Lyapunov quantities using the algorithm of Section~\ref{se:lyapunovconstants}. Recall that, as we have a $CC$-equilibrium point, $W_1 = 0.$ Then, the necessary conditions on the parameters to have a center follow from solving (at least) the algebraic system with four equations and three parameters defined by
\begin{equation*}
	\mathcal{S} = {W_2 = W_4 = W_6 = W_8 = 0}.
\end{equation*}
All expressions are, by definition, considered when the previous quantities vanish. So, the values of $W_3,W_5,$ and $W_7$ vanish.

Here, instead of presenting the complete expressions of $W_i$, $i=2,4,6,8$, for simplicity, we provide only some numerators due to the size of these expressions, denoted by $\widehat W_i$, which are polynomials with integer coefficients in the remaining parameters $(e_1,b_2,e_2)$. In particular, from the last one, the different families in the statement can be seen by analyzing when each factor of $\widehat W_8$ vanishes in the reals. The numerators of $\widehat W_4$ and $\widehat W_6$ are polynomials of degree $20$ and $104$, with $295$ and $1832$ monomials, respectively. Hence, we have:
\begin{equation*}
\begin{aligned}
\widehat W_2=&
27 b_2 e_1 (9 e_1^2 - 24 e_1 + 32) e_2^4 - (1215 b_2^2 e_1^3 - 3240 b_2^2 e_1^2 + 243 e_1^4 + 4320 b_2^2 e_1 \\
& - 1215 e_1^3 + 2916 e_1^2 - 3456 e_1 + 2304) e_2^3 + 12 b_2 (189 b_2^2 e_1^3 - 504 b_2^2 e_1^2 + 54 e_1^4 \\
& + 672 b_2^2 e_1 - 216 e_1^3 + 504 e_1^2 - 576 e_1 + 512) e_2^2 - 16 (b_2^2 + 1) (117 b_2^2 e_1^3 \\
& - 312 b_2^2 e_1^2 + 27 e_1^4 + 416 b_2^2 e_1 - 144 e_1^3 + 348 e_1^2 - 416 e_1 + 256) e_2 \\
&+ 64 e_1 b_2 (9 e_1^2 + 24 e_1 + 32) (b_2^2 + 1)^2,\\
\widehat W_8=&e_1^5(1458e_1^8 - 15309e_1^7 + 78003e_1^6 - 248832e_1^5 + 535680e_1^4- 781056e_1^3  + 755712e_1^2\\ 
&- 442368e_1 + 131072)^3(b_2+1)(3e_1 - 8)^2(3e_1^2 - 7e_1 + 8)^2(27e_1^4 - 108e_1^3   \\
	&+ 288e_1^2- 384e_1+ 256)b_2^3(b_2^2 + 1)^7(9e_1^2 - 24e_1 + 32)^4(3e_1 - 4)^2(b_2- 1)\\
	&(7015420248855589903862275823952189119801685375120618e_1^{11} \\
	&- 152357242285056016622382435380491647911148728341085013e_1^{10}\\
	& + 1435336366582110818156465089095809417442072336316281705e_1^9\\
	& - 8043276072604152065018524298720295090916067001323121132e_1^8 \\
	&+ 30283544465057925411860791413515197107177223647053806704e_1^7\\
	& - 81061414503281949781688061029238613292150397149292296448e_1^6 \\
	&+ 158934347738141481809466504567893193833879482270751697792e_1^5 \\
	&- 230418094422134302482528539805203263325388508317820726272e_1^4 \\
	&+ 243754882522350777561809342871138194082313826638503997440e_1^3 \\
	&- 179413407752255527803949903705145705530859993913011601408e_1^2\\
	&+ 81749048336969359284906573839660434291536221323282087936e_1\\
	&-16511335007295362432559652316381263201043632657029660672).
\end{aligned}
\end{equation*}
The proof finishes by checking that $C_1, \ldots, C_8$ are also sufficient conditions for having a center, using that $Z_i$, $i = 1,2$, have a first integral of the form \eqref{ip}. Under the conditions $C_1$ and $C_2$, the first integrals coincide on the separation line, so it follows from Proposition~\ref{teocentros} that the monodromic equilibrium point is a center. For the remaining conditions, it is necessary to manipulate it before applying again Proposition~\ref{teocentros}. For case $C_3$, for example, the expressions of the first integrals in each zone can be written, up to a multiplicative constant, as
\[
\begin{aligned}
	H_1(x,y) &= \frac{(3y)^{1/4}(20x - 3y - 5)}{4x^{5/3}}, \\
	H_2(x,y) &= \frac{(3y)^{(b_2^2 - 3)/4}\left(4b_2^2x + (3 - b_2^2)(3y + 1)\right)}{12x^{b_2^2/3}},
\end{aligned}
\]
for $i = 1, 2$. Defining $\widetilde{H}_1 = H_1^{12/5}$ and $\widetilde{H}_2 = H_2^{12/(b_2^2 - 3)}$, both first integrals coincide on the straight line of separation:
\[
\widetilde{H}_i\left(x, \dfrac{4(x - 1)}{3} + 1\right) = \frac{(4x - 1)^3}{x^4}.
\]
Consequently, we obtain a continuous first integral. The other cases follow analogously.
\end{proof}
	
\section{Bifurcation of limit cycles in $\mathcal K_2$} \label{bifurcation}
In this section, we study the unfolding of small amplitude crossing limit cycles from a weak focus on $\Sigma.$ We consider the perturbation of a $CC$-equilibrium point with two different approaches, but in both cases we perturb an analytic quadratic Kolmogorov center. Firstly, we perturb maintaining (partially) integrable systems on both sides of the separation straight line. Secondly, we break this integrability condition. All limit cycles are of small size and bifurcate from a $CC$-equilibrium point. In the first case, we have obtained weak foci of eighth-order, but due to the special structure of the perturbation family, only four or five limit cycles bifurcate depending on the perturbation class. In the second case, we have obtained a higher number of limit cycles, six, from weak foci of sixth-order that unfold in a versal way. Of course, under generic perturbation, we could unfold eight limit cycles, but the restriction to being in the Kolmogorov class makes the bifurcation analysis more intricate. We recall that, in the piecewise scenario, we have previously explained that from a weak focus of order $k,$ one can generically unfold $k$ limit cycles using an analytic perturbation. See \cite{GouvTorre2021l} for more details.

\subsection{Bifurcation of limit cycles in the center-center case}
This section is devoted to the study of the limit cycles of crossing type in the class of piecewise quadratic Kolmogorov systems, $\mathcal{K}_2,$ but restricting the analysis to systems having a center in each piece. In Proposition~\ref{bclcc}, we analyze the bifurcation problem of limit cycles from a $CC$-equilibrium point in this class but without sliding segment. This restriction is partially not considered in Corollary~\ref{co1}, getting one extra limit cycle. The bifurcation mechanism is the one introduced in Section~\ref{se:Lyapunovconstantsseries} from the computation of the first-order terms of the Lyapunov quantities.

\begin{pro}\label{bclcc} 	
Consider the piecewise perturbed polynomial system of the form \eqref{eq:centrob}, given by:
\begin{equation}\label{eq:ps}
\begin{aligned}
	X_{c} &= x \left(bx - \frac{(b^2 + 1)}{e} y + \frac{b^2 - be + 1}{e}\right), \\	
	Y_{c} &= y \left(ex - by + b - e\right), \\
	X_i &= x \left(q_i x + \frac{(1+b^2)p_i - (2b + q_i)eq_i }{(e + p_i)e} y - \frac{(b^2 + e q_i + 1)p_i + (e - 2b   -  q_i)eq_i}{(e + p_i)e}\right), \\
	Y_i &= y \left(p_i x - q_i y + q_i - p_i\right),
\end{aligned}
\end{equation}
where \((e + p_i)e \neq 0\) for \(i = 1, 2\), \(\mu = (b, e) \in \mathbb{R}^2\), \(\lambda = (p_1, q_1, p_2, q_2) \in \mathbb{R}^4\) are small enough, and the separation straight line is \( 4(x - 1) - 3(y - 1)=0\). Then, there exist values of the parameters for this perturbed polynomial system that yield a eighth-order weak focus at \((1, 1)\), which can unfold at least three crossing limit cycles of small amplitude.
\end{pro}	

\begin{proof}
System~\eqref{eq:ps}, with \(p_1 = q_1 = p_2 = q_2 = 0\), is the same on both sides of the separation line. Moreover, it has a center at \((1,1)\) because it possesses a first integral of the form \eqref{ip}, satisfying that the trace is zero and the determinant is one. Consequently, in this case, we have a center for all separation lines passing through the \(CC\)-equilibrium point. Each piece defining \eqref{eq:ps}, for \(i=1,2\), is integrable of center type, but the center conditions for the piecewise system depend on the values of the parameters \(\lambda\).
	
By construction, system \eqref{eq:ps} has no sliding segment, so \(W_0=W_1 = 0\) and \((1,1)\) is an equilibrium point of weak focus type. Its stability can be determined using the algorithm of Section~\ref{se:lyapunovconstants}, but first, we need to perform a translation and a rotation so that the equilibrium point is located at the origin and the separation line \(4x - 3y = 0\) becomes \(y = 0\). Straightforward computations provide the first Lyapunov quantities \(W_j\) for \(j = 2, \ldots, 8\). As shown in Section~\ref{se:lyapunovconstants}, the Lyapunov quantities with odd subscripts are zero when the previous also are; only those with even subscripts can be used to study the number of zeros of the difference map \eqref{eq:difmap}. We find that they are rational functions that admit Taylor expansions in \(\lambda\), and the first-order jet can be computed, yielding, in this case:
\begin{equation*}
	\begin{aligned}
		W^{[1]}_{2} =& \frac{2e}{243} L_{2}(b,e) (q_1 - q_2) + \frac{2}{243} M_{2}(b,e) (p_1 - p_2),
	\end{aligned}
\end{equation*}
where
\begin{equation*}
	\begin{aligned}
		L_{2} &= -243 e^6 - 1458 b e^5 + 8(27 - 31 b^2) e^4 + 288 (4 b^2 + 3)b e^3+ 192(42 b^4 + 27 b^2 - 5) e^2  \\
		& \quad + 512b(17 b^2 + 5)(b^2 + 1) e + 1024 (3b^2 + 1) (b^2 + 1)^2, \\
		M_{2}&= -243 b e^6 - 1296 b^2 e^5 - 972 b (3 b^2 + 2) e^4 - 2016 (2 b^2 + 1) (b^2 + 1) e^3\\
		& \quad  - 192 b (22 b^2 + 17) (b^2 + 1) e^2 - 3072 b^2 (b^2 + 1)^2 e - 1024 b (b^2 + 1)^3.
	\end{aligned}
\end{equation*}
Assuming that \(L_2(b,e) \ne 0\) and using the Implicit Function Theorem, we can express \(q_1\) in terms of $(p_1, p_2, q_2)$. Therefore, the first-order jet of the next Lyapunov quantities, for \(l = 2, 3, 4\), can be computed as
\begin{equation*}
W^{[1]}_{2l} = \frac{M_{2l}(b,e)}{L_{2}(b,e)} (p_1 - p_2),
\end{equation*}
where
\begin{equation*}
M_{2l}(b,e) = e^{2l} (4b + 3e)^2 (4b^2 - be - 3e^2 + 4)^2 (4b^2 + 3be + 4)^2 (16b^2 + 24be + 9e^2 + 16)^{2l} m_{2l}(b,e),
\end{equation*}
with
\begin{equation*}
\begin{aligned}
	m_{4}(b,e) &= -27 b e^4 - (81 b^2 - 162) e^3 - 36 b^3 e^2 + (8 b^4 - 192 b^2 - 272) e  + 64 b (b^2 + 1)^2, \\
	m_{6}(b,e) &= -346428 b^3 e^{12} - (31177872 b^2 - 149328 ) b^2e^{11} - (1739448 b^4 - 8485776 b^2 \\
	& \quad - 36551331 ) be^{10} - (145881648 b^6 - 172851732 b^4 - 9128457 b^2 - 33736662) e^9 \\
	& \quad + (6624864 b^6 + 1459224 b^4 - 2776284 b^2 - 162724464 ) be^8 + (451638 b^8 \\
	& \quad + 1337868 b^6 - 12979584 b^4 - 88322148 b^2 - 223416) e^7 + (434294784 b^9 \\
	& \quad + 296182656 b^6 - 46674 b^4 - 2553264 b^2 + 1383552) be^6 - (167878656 b^{10} \\
	& \quad + 86413824 b^8 - 1349768448 b^6 - 26465832 b^4 - 1792233216 b^2 - 414288) e^5 \\
	& \quad - 27648 (b^2 + 1)^2 (2768 b^6 - 2678 b^4 - 3559 b^2 - 15344) be^4 - 496 (68288 b^6 \\
	& \quad + 36774 b^4 + 73641 b^2 + 58141) (b^2 + 1)^3 e^3 + 49152 (2816 b^4 - 5634 b^2 \\
	& \quad - 6551) (b^2 + 1)^4 b e^2 + 393216 (b^2 + 1)^5 (44 b^2 - 279) b^2 e + 46137344 b^3 (b^2 + 1)^6,
\end{aligned}
\end{equation*}
and \(m_{8}(b,e)\) is a polynomial of degree 25. We have not presented it explicitly due to its size. The solutions of the algebraic system \(\{m_{4}(b,e) = 0, m_{6}(b,e) = 0\}\) can be explicitly written as
\begin{equation}\label{pt1}
	(b^*, e^*) = \left(\gamma, - \frac{\gamma(1392580 \gamma^4 + 544194 \gamma^2 - 1145761)}{62169}\right),
\end{equation}
where \(\gamma\) is each real root of the polynomial
\begin{equation*}
	g(\gamma) = 9604 \gamma^6 - 1470 \gamma^4 - 9797 \gamma^2 + 4232.
\end{equation*}
Clearly, \(W^{[1]}_{4}(b^*, e^*) = W^{[1]}_{6}(b^*, e^*) = 0\), and then, we can check that the next coefficient is nonvanishing because its resultant with \(g\) with respect to \(\gamma\) is a non-zero rational number. In fact, we can obtain
\begin{equation*}
\begin{aligned}
	W^{[1]}_{8}(b^*, e^*) &= \frac{753747879498059507302400000}{10557} \gamma^5 + \frac{14617417205145737048883200000}{517293} \gamma^3 \\
	& \quad - \frac{3294809732513451317657600000}{57477} \gamma.
\end{aligned}
\end{equation*}
Hence, we have curves of weak foci of eighth-order passing through each point of the form \eqref{pt1}. Moreover, the intersection of \(m_4\) and \(m_6\) is transversal at such points because the determinant of the Jacobian matrix is
\begin{equation*}
\begin{aligned}
	\det\Jac(b^*, e^*) &= \frac{1163618760932343152640}{901} \gamma^4 + \frac{3223719345310777999360}{6307} \gamma^2 \\
	& \quad - \frac{6539723921076330168320}{6307}.
	\end{aligned}
\end{equation*}
Again, the nonvanishing condition of the above polynomial is guaranteed after checking that its resultant with \(g\), with respect to \(\gamma\), is a non-zero rational number. So, using Proposition~\ref{prop:theop}, the proof follows.
\end{proof}

As we have commented, the above result has no sliding segment. With Proposition~\ref{hopf_bif}, we can add an extra limit cycle of small amplitude using a Hopf-type bifurcation while remaining within the $\mathcal{K}_2$ class. Hence, after an adequate perturbation, we will have a sliding or escaping segment. This description proves the next immediate consequence of Proposition~\ref{bclcc}.
\begin{coro}\label{co1}
There exist $Z\in\mathcal K_2,$ with separation straight line $4(x-1)-3(y-1)=0$ and a weak focus of eighth-order at $(1,1),$ which unfold at least four crossing limit cycles of small amplitude in the center-center case.
\end{coro}

The previous result naturally motivates us to investigate the focus-focus case to achieve a complete unfolding, obtaining eight limit cycles in total. The proofs of the above two results make it clear that only half of the coefficients of the difference map play a role during the unfolding procedure to obtain a degenerate Hopf bifurcation. This is similar to what occurs in the analytic case. If we move beyond the center-center case by considering nonvanishing traces, the computations used to obtain the coefficients $W_k$
are not valid. However, we will discuss this problem extensively in the next section. From the proof of Proposition~\ref{bclcc}, it is evident that by adding only the trace parameters, a fifth limit cycle of small amplitude could be obtained in Corollary~\ref{co1}. But as more limit cycles can appear in the next section, we have not added this particular result here.

\subsection{Bifurcation of limit cycles in the focus-focus case}
In this section, we will focus on studying limit cycles of the crossing type when we do not have first integrals in each piece defined by systems in the $\mathcal{K}_2$ class. That is, when both systems have equilibria (visible or invisible) of focus type, we refer to this as the focus-focus case. We will follow the same scheme as in the previous section, starting with a result, Proposition~\ref{promean}, that does not involve a sliding or escaping segment, and although the order of the weak focus is lower, the unfolding is versal, leading to higher cyclicity by perturbing the same center as in the previous section. The computations are more intricate because we need to analyze the jets of order two. As a consequence, using Proposition~\ref{hopf_bif}, the proof of our main result, Theorem~\ref{th:m1}, follows by again adding a limit cycle due to the birth of an escaping or sliding segment. Note that for the following result, after further investigation, we noticed that introducing additional perturbative monomials does not result in more limit cycles.
	
\begin{pro}\label{promean} Consider the piecewise polynomial perturbed system of the form \eqref{eq:a8}, with 
\begin{equation}\label{eq:9}
	\begin{aligned}
		(X_{c},Y_{c})=&(x(bx-(b^2 + 1)y/e +( b^2-be+1)/e),y(ex-by+b-e)),\\
		(X_i,Y_i)=&(p_{i0}x(x - 1 ),q_{i1}y( y - 1)),
	\end{aligned}
\end{equation}	
where  $e\neq0,$  $\mu=(b,e)\in\R^2,$ $\lambda=(p_{10},p_{20},q_{11},q_{21})\in\R^4,$ with $i=1,2$ small enough, and the straight separation line is $4(x-1)-3(y-1)=0.$ Then, there exist values of the parameters such that the perturbed system has a curve of weak foci of sixth-order at $(1,1),$ which can unfold at least five crossing limit cycles of small amplitude.
\end{pro}		

\begin{proof}
The perturbed system in the statement has a monodromic $CC$-equilibrium at $(1,1)$. In order to use the mechanism described in Section~\ref{se:Lyapunovconstantsseries}, we translate $(1,1)$ to the origin and subsequently apply the affine transformation to write the Jacobian matrix at the equilibrium in its Jordan normal form. Then, the expressions in \eqref{eq:9} become
\begin{equation*}
\begin{aligned}
	X_{c}=& y+bx^2 + (b^3 - b^2e + b + e)yx/e + by^2,\\
	Y_{c}=&-x -(b^2 + 1)xy/e,\\
	X_i=&q_{i1}x- (p_{i0} - q_{i1})b y + 2q_{i1}bxy- (p_{i0} b^2 - q_{i1}b + p_{i0})by^2/e+q_{i1}x^2,\\
    Y_i=&p_{i0}y(1  + (b^2  + 1)y/e).
\end{aligned}
\end{equation*}	
Now, the origin of the unperturbed system is a non-degenerate center. We will need to carry out a second-order analysis of the first coefficients of the difference map, proving that their Taylor expansions (of order two), denoted by $\mathcal{W}^{[2]}_{j},$  are transversally independent. Hence,  we can write the difference function as
\begin{equation*}
\Delta(r)=\sum_{j=1}^6w_jr^j+O(r^7),
\end{equation*} 
where $w_j\in\R,$ for $j=1,\ldots,6.$ Straightforward computations show that some of the linear parts are linearly independent, allowing us to perform a change of variables in the parameter space to adopt them as new parameters. We then write
\begin{equation}\label{secondorder}
	\begin{aligned} 
	\mathcal{W}^{[2]}_{1}=&\omega_1+W^{[2]}_{1}, &
	\mathcal{W}^{[2]}_{2}=&\omega_2+W^{[2]}_{2}, &
	\mathcal{W}^{[2]}_{3}=&W^{[1]}_{3}+W^{[2]}_{3},\\
	\mathcal{W}^{[2]}_{4}=&\omega_4+W^{[2]}_{4}, &
 	\mathcal{W}^{[2]}_{5}=&W^{[1]}_{5}+W^{[2]}_{5},& \mathcal{W}^{[2]}_{6}=&\omega_6+W^{[2]}_{6}, \\
	\end{aligned}
\end{equation}
with 
\begin{equation}\label{seconorder1}
\begin{aligned} 
W^{[1]}_{3}=&\frac{1}{2^43^2 5^6}\frac{L_3(b,e)}{e^2}\omega_1-\frac{1}{3(5)^3}\frac{M_3(b,e)}{e}\omega_2, \\
W^{[1]}_{5}=&\frac{1}{2^83^45^{12}}\frac{L_5(b,e)}{e^4}\omega_1-\frac{1}{3^35^9}\frac{M_5(b,e)}{e^3}\omega_2
-\frac{2}{5^3}\frac{R_5(b,e)}{e}\omega_4,
\end{aligned}
\end{equation}
where $L_3,M_3, L_5, M_5,$ and $R_5,$ have degree 6, 3, 12, 9, and 3, respectively. Moreover, each $W^{[2]}_{j}$ is a rational function $W^{[2]}_{j}=N^{[2]}_{j}/D^{[2]}_{j}(b,e)$  where $N^{[2]}_{j}$ are, by construction, homogeneous polynomials of degree $2$ in $\omega_1,$ $\omega_2,$ $\omega_4,$  $\omega_6,$ with integer coefficients. The denominators $D^{[2]}_{j}(b,e)$ are polynomials in $b, e,$ with shared factors but different exponents. Due to their size, we omit them here. Additionally, the next two coefficients $W^{[1]}_{7}, W^{[1]}_{8}$ are expressed as $W^{[1]}_{5},$ i.e., linear functions in $\omega_1, \omega_2,$ and $\omega_4.$ Then, using first-order analysis via Proposition~\ref{thm:chrorder1}, we obtain $3$ small-amplitude limit cycles. 

The next step, equivalent to applying the Implicit Function Theorem, is to simplify \eqref{seconorder1} using \eqref{secondorder}, eliminating their linear parts by affine transformations. The resulting second-order terms are:
\begin{equation*}
	\begin{aligned}
		\widetilde{\mathcal{W}}^{[2]}_{3}=&\mathcal{W}^{[2]}_{3}-\left(\frac{1}{2^43^25^6}\frac{L_3(b,e)}{e^2} \mathcal{W}^{[2]}_{1}  -\frac{1}{3(5)^ 3}\frac{M_3(b,e)}{e} \mathcal{W}^{[2]}_{2}\right),\\
		\widetilde{\mathcal{W}}^{[2]}_{5}=&\mathcal{W}^{[2]}_{5}-\left(\frac{1}{2^8 3^45^{12}}\frac{L_5(b,e)}{e^4}\mathcal{W}^{[2]}_{1}-\frac{1}{3^35^9}\frac{M_5(b,e)}{e^3}\mathcal{W}^{[2]}_{2}-\frac{2}{5^3}\frac{R_5(b,e)}{e}\mathcal{W}^{[2]}_{4}\right),
	\end{aligned}
\end{equation*} 
where 
$\widetilde{\mathcal{W}}^{[2]}_{j}=\widetilde{N}^{[2]}_{j}/\widetilde{D}^{[2]}_{j}(b,e)$. Again, $\widetilde{N}^{[2]}_{j}$ are homogeneous polynomials of degree $2$ in $\omega_1,$ $\omega_2,$ $\omega_4,$ and $\omega_6,$ with rational coefficients, and $\widetilde{D}^{[2]}_{j}(b,e)$ are polynomials in $b,e,$ for $j=3,5$. The transversality of the second-order terms is not affected by the denominators, which are nonvanishing. Then \eqref{secondorder} becomes  
\begin{equation}\label{cv}
	\begin{aligned}
		\mathcal{W}^{[2]}_{1}=&\omega_1+\frac{N^{[2]}_{1}}{D(b,e)}, \quad \mathcal{W}^{[2]}_{2}=\omega_2+\frac{N^{[2]}_{2}}{D(b,e)}, \quad 
		\widetilde{\mathcal{W}}^{[2]}_{3}=\frac{\widetilde{N}^{[2]}_{3}}{D(b,e)}, \\
		 \mathcal{W}^{[2]}_{4}=&\omega_4+\frac{N^{[2]}_{4}}{D(b,e)}, \quad 
		 \widetilde{\mathcal{W}}^{[2]}_{5}=\frac{\widetilde{N}^{[2]}_{5}}{D(b,e)}, \quad \mathcal{W}^{[2]}_{6}=\omega_6+\frac{N^{[2]}_{6}}{D(b,e)},
	\end{aligned}
\end{equation} 
where $D(b,e),$ is the nonvanishing common denominator. Now, consider new weighted variables, $\omega_1=\omega^2 \widetilde{\omega}_1,$  $\omega_2=\omega^2 \widetilde{\omega}_2,$  and $\omega_4=\omega^2 \widetilde{\omega}_4,$ $\omega_6=\omega.$ Substituting into \eqref{cv}, we rewrite it as 
\begin{equation*}
	\begin{aligned}
		\mathcal{W}^{[2]}_{1}=&\omega^2\widetilde{\omega}_1+O(\omega^3), & \mathcal{W}^{[2]}_{2}=&\omega^2 \widetilde{\omega}_2+O(\omega^3),\\
		\widetilde{\mathcal{W}}^{[2]}_{3}=&\omega^2 e^8 R(b,e)\frac{L(b,e)}{D(b,e)}+O(\omega^3),&
		\mathcal{W}^{[2]}_{4}=&\omega^2 \widetilde{\omega}_4+O(\omega^3),\\
		 \widetilde{\mathcal{W}}^{[2]}_{5}=&\omega^2e^6R(b,e)\frac{M(b,e)}{D(b,e)}+O(\omega^3),&
		 \mathcal{W}^{[2]}_{6}=&\omega+O(\omega^2),
	\end{aligned}
\end{equation*} 
with $$R(b,e)=63 b^3-63 b^2 e-4 b^2+8 b e+63 b+18 e-4,$$
where $L(b,e)$ and $M(b,e)$ are polynomials of degree 18 and 24, respectively. We find 14 solutions of the system $\{L(b,e)=M(b,e)=0\}$ under the condition $D(b,e)\neq0,$ and only one yields transversal intersection points for each root of a given degree 36 polynomial $f.$ In fact, there exist a polynomial $g$ such that the intersection is written as $(b^*,e^*)=(z,g(z))$ being $z$ a simple zero of $f.$  One of these solutions is $(b^*,e^*)\approx(-0.3589344145,1.09217769345).$ The transversality is verified by  checking that the determinant of the Jacobian matrix and the denominator $D$ do not vanish at this point. These last properties are confirmed  by computing the resultants with $f$ with respect to $z$ and verifying that they are nonzero. The existence of a transversal curve of sixth-order weak foci follows by taking $(b^*,e^*),$ and $\widetilde{\omega}_1=\widetilde{\omega}_2=\widetilde{\omega}_4=0,$ for small enough $\omega$. The proof follows again using Proposition~\ref{prop:theop}.
\end{proof}
	
Other centers in Proposition~\ref{center2} have also been considered and their bifurcation unfolding studied within Kolmogorov quadratic family, but fewer small-amplitude limit cycles were obtained. Therefore, we have not detailed them explicitly here.

\section{Conclusions and further remarks}	
In this work, we study the classical center-focus and cyclicity problems in the setting of piecewise quadratic Kolmogorov systems. For the first problem, we assume a natural condition on the equilibrium: it lies on the separation line and is of CC-equilibrium type. We provide a partial classification of centers for different separation lines. Achieving a complete classification, however, requires a considerable computational effort, in contrast with the smooth case where the solution is straightforward.

For the second problem, unlike the smooth quadratic case where no limit cycles exist, the piecewise scenario allows for the existence of at least six crossing limit cycles. Moreover, we identify foci of high weak order, eight in fact, that do not versally unfold into this maximum number of bifurcated limit cycles within the Kolmogorov class. This illustrates both the computational challenges and the distinctive differences between the two problems: the number of small-amplitude limit cycles and the order of weak foci in the piecewise context.

In Section~\ref{bifurcation}, we restricted our analysis to the particular separation straight line $\Sigma$ taking $(\alpha,\beta)=(4,3)$. We have obtained similar results for other choices, $(\alpha,\beta)\in\{(1,2), (2,1), (3,1)\}$, which we have not included in order to avoid unnecessary repetitions. These cases reinforce the intuition that higher-order weak foci or additional small-amplitude limit cycles will not appear. We remark that the general problem is extremely intricate, but we do not expect it to yield new or qualitatively different results.

These findings highlight the richer dynamics of non-smooth systems. In particular, it has been necessary to adapt previous results to analyze the conditions that guarantee the existence of a pseudo-Hopf bifurcation in the framework of piecewise Kolmogorov systems. We believe that the technical tools and algorithms developed here will be useful for future studies on bifurcation phenomena, especially in contexts where classical bifurcation techniques must be adapted.
	
\section*{Acknowledgements}
This work has been realized thanks to the funded by the Brazilian S\~ao Paulo Research Foundation FAPESP grants  2021/14987-8, 2022/14484-9 and, 2021/21181-0;  the National Council for Scientific and Technological Development(CNPq) grants 407454/2023 and 304766/2019-4; the Catalan AGAUR Agency (grant 2021 SGR 00113); the Spanish AEI agency (grants PID2022-136613NB-I00 and CEX2020-001084-M).

\bibliographystyle{abbrv}
\bibliography{biblio} 
\end{document}